\documentclass[11pt]{article}

\usepackage{amssymb}

\makeatletter\@addtoreset{equation}{section}

\newtheorem{theorem}{Theorem}[section]
\newtheorem{corollary}[theorem]{Corollary}
\newtheorem{lemma}[theorem]{Lemma}
\newtheorem{proposition}[theorem]{Proposition}

\def\endverif{\nopagebreak\newline\mbox{\ }\hfill\rule{2mm}{2mm}}
\def\endveriff{\nopagebreak\mbox{\ }\hfill\rule{2mm}{2mm}}

\def\C{{\mathbb C}}
\def\N{{\mathbb N}}
\def\T{{\mathbb T}}
\def\Z{{\mathbb Z}}

\def\NN{{\cal N}}
\def\RT{{\cal R}_T}

\def\Ad{{\rm Ad}\,}
\def\Aut{{\rm Aut}}
\def\AutT{{\rm Aut}(R_T,S_T\,|\,C_T)}
\def\AutTT{{\rm Aut}(R_T,C_T\,|\,S_T)}
\def\AutX{{\rm Aut}(X,\mu)}
\def\id{{\rm id}}
\def\Int{{\rm Int}}
\def\RR{{\rm Re}\,}

\def\<{\langle}
\def\>{\rangle}
\def\b{\overline}

\def\eps{\varepsilon}

\begin{document}

\title{\bf Ergodic theory and maximal abelian subalgebras of the
hyperfinite factor}

\author{Sergey Neshveyev$^{1)}$ and Erling St{\o}rmer$^{2)}$}

\date{}

\footnotetext[1]{Partially supported by the Royal Society/NATO
postdoctoral fellowship, the Centre for Advanced Study in Oslo 
and Award No UM1-2092 of the Civilian
Research \& Development Foundation}

\footnotetext[2]{Partially supported by the Norwegian
Research Council}

\maketitle

\begin{abstract}
Let $T$ be a free ergodic measure-preserving action of an abelian
group $G$ on $(X,\mu)$. The crossed product algebra
$R_T=L^\infty(X,\mu)\rtimes G$ has two distinguished masas, the
image $C_T$ of $L^\infty(X,\mu)$ and the algebra $S_T$ generated
by the image of $G$. We conjecture that conjugacy of the singular
masas $S_{T^{(1)}}$ and $S_{T^{(2)}}$ for weakly mixing actions
$T^{(1)}$ and $T^{(2)}$ of different groups implies that the
groups are isomorphic and the actions are conjugate with respect
to this isomorphism. Our main result supporting this conjecture is
that the conclusion is true under the additional assumption that
the isomorphism $\gamma\colon R_{T^{(1)}}\to R_{T^{(2)}}$ such
that $\gamma(S_{T^{(1)}})=S_{T^{(2)}}$ has the property that the
Cartan subalgebras $\gamma(C_{T^{(1)}})$ and $C_{T^{(2)}}$ of
$R_{T^{(2)}}$ are inner conjugate. We discuss a stronger
conjecture about the structure of the automorphism group
$\Aut(R_T,S_T)$, and a weaker one about entropy as a conjugacy
invariant. We study also the Pukanszky and some related invariants
of $S_T$, and show that they have a simple interpretation in terms
of the spectral theory of the action $T$. It follows that
essentially all values of the Pukanszky invariant are realized by
the masas $S_T$, and there exist non-conjugate singular masas with
the same Pukanszky invariant.
\end{abstract}

\section{Introduction}
It is well-known that if one has a Lebesgue space $(X,\mu)$ with a
free ergodic measure-preserving action $T$ of an abelian group
$G$, then the crossed product algebra $R_T=L^\infty(X,\mu)\rtimes
G$ is the hyperfinite factor with two distinguished maximal
abelian subalgebras (masas), the image $C_T$ of $L^\infty(X,\mu)$
and the masa $S_T$ generated by the canonical unitaries in $R_T$
implementing the action. It is the purpose of the present work to
investigate how much information about the system $(X,\mu,T)$ can
be extracted from properties of the masas $C_T$ and $S_T$.

In Section~\ref{2} we formulate our main conjecture that for
weakly mixing actions the masas $S_T$ determine the actions up to
an isomorphism of the groups. Here we also give a short proof of
the singularity of $S_T$, a result due to Nielsen~\cite{Nielsen},
and more generally describe the normalizer of $S_T$ for arbitrary
actions, which is a result of Packer~\cite{Packer1}.

Apparently the only conjugacy invariant of singular masas which
has been effectively used over the years, is the invariant of
Pukanszky~\cite{Pukanszky}. It arises as a spectral invariant of
two commuting representations of a masa $A\subset M$ on
$B(L^2(M))$ coming from the left and right actions of $A$ on $M$.
It is not surprising that for the masas $S_T$ this invariant is
closely related to spectral properties of the action $T$. This
fact has two consequences. On one hand, we have a lot of actions
with different Pukanszky invariants. On the other hand, for most
interesting systems such as Bernoullian systems, the invariant
gives us nothing. This is described in Section~\ref{3}.

In Section~\ref{4} we prove the main result supporting our
conjecture. Namely, for weakly mixing actions the pair consisting
of the masa $S_T$ and the inner conjugacy class of $C_T$ is an
invariant of the action. In fact if $\Aut(R_T,S_T)$ denotes the
subgroup of $\gamma\in\Aut(R_T)$ such that $\gamma(S_T)=S_T$, we
prove a stronger result describing the subgroup of $\Aut(R_T,
S_T)$ consisting of automorphisms~$\gamma$ such that $\gamma(C_T)$
and $C_T$ are inner conjugate. We conjecture that this subgroup is
actually the whole group $\Aut(R_T, S_T)$. One test for our
conjecture is to prove that this subgroup is closed, and we are
able to do this under slightly stronger assumptions than weak
mixing.

The group of inner automorphisms defined by unitaries in $S_T$ is
not always closed, and this gives us the possibility of
constructing non-conjugate singular masas with the same Pukanszky
invariant.

Finally in Section~\ref{5}, which is independent of the others, we
consider a weaker conjecture stating that the entropy of the
action is a conjugacy invariant for $S_T$. We prove that if
$S_{T^{(1)}}$ and $S_{T^{(2)}}$ are conjugate and under this
conjugacy the canonical generators of these algebras coincide on a
small projection, then the entropies of the actions coincide. The
proof is an application of the theory of non-commutative entropy.

\medskip\noindent
{\it Acknowledgement.} The authors are indebted to J.~Packer,
S.~Popa and A.~Vershik for helpful discussions.

\bigskip\bigskip

\section{Preliminaries on crossed products} \label{2}

Let $G$ be a countable abelian group, $g\mapsto T_g\in\Aut(X,\mu)$
a free ergodic measure-preserving action of $G$ on a Lebesgue
space $(X,\mu)$. Consider the corresponding action
$g\mapsto\alpha_g$ on $L^\infty(X,\mu)$, $\alpha_g(f)=f\circ
T_{-g}$, and the crossed product algebra
$L^\infty(X,\mu)\rtimes_\alpha G$, which will be denoted by~$R_T$
throughout the paper. Let $g\mapsto v_g$ be the canonical
homomorphism of $G$ into the unitary group of $R_T$. We denote by
$S_T$ the abelian subalgebra of $R_T$ generated by $v_g$, $g\in
G$. The algebra $L^\infty(X,\mu)$ considered as a subalgebra of
$R_T$ will be denoted by $C_T$.

To fix notations, the unitary on $L^2(Y,\nu)$ associated with
an invertible non-singular transformation $S$ of a measure space
$(Y,\nu)$ will be denoted by $u_S$,
$u_Sf=(dS_*\nu/d\nu)^{1/2}f\circ S^{-1}$, and the corresponding
automorphism of $L^\infty(Y,\nu)$ will be denoted by $\alpha_S$,
$\alpha_S(f)=f\circ S^{-1}$. For a given action $T$ we shall
usually suppress $T$ in such notations, so we write $u_g$ and
$\alpha_g$ instead of $u_{T_g}$ and $\alpha_{T_g}$.

We shall usually consider $R_T$ in its standard representation on
$L^2(X,\mu)\otimes L^2(\hat G,\lambda)$, where $\hat G$ is the
dual group and $\lambda$ is its Haar measure. The elements of the
group $G$ considered as functions on $\hat G$ define two types of
operators on $L^2(\hat G)$, the operator $m_g$ of multiplication
by~$g$, $(m_gf)(\chi)=\<\chi,g\>f(\chi)$, and the projection $e_g$
onto the one-dimensional space $\C g$. Then the representation
$\pi$ of $R_T$ on $L^2(X)\otimes L^2(\hat G)$ is given by
$$
\pi(v_g)=1\otimes m_g,\ \ \pi(f)=\sum_g\alpha_g(f)\otimes e_{-g}\
\hbox{for}\ f\in L^\infty(X).
$$
Then $R_T$ is in its standard form with the tracial vector
$\xi\equiv1$. The modular involution $J$ is given by
\begin{equation} \label{e2.1}
J=\tilde J\sum_g u_g\otimes e_g=\left(\sum_g u_g\otimes
e_{-g}\right)\tilde J,
\end{equation}
where $\tilde J$ is the usual complex conjugation on
$L^2(X\times\hat G)$. Indeed, since $(1\otimes e_g)\xi=0$ for
$g\ne0$,
\begin{eqnarray*}
J\pi(v_gf)\xi
 &=&J\pi(v_g)(f\otimes e_0)\xi=J(f\otimes m_ge_0)\xi=J(f\otimes e_gm_g)\xi
      =\tilde J(u_gf\otimes e_gm_g)\xi\\
 &=&\tilde J(\alpha_g(f)\otimes e_gm_g)\xi=(\alpha_g(\b{f})\otimes
      e_{-g}m_{-g})\xi=\pi(\b{f}v_g^*)\xi.
\end{eqnarray*}
In particular,
\begin{equation} \label{e3.1}
J\pi(f)J=\overline{f}\otimes1,\ \ J\pi(v_g)J=u_g\otimes
m_g^*.
\end{equation}
Hence $R_T$ is the fixed point subalgebra of $L^\infty(X)\otimes
B(L^2(\hat G))$ for the action $g\mapsto\alpha_g\otimes\Ad m_g^*$
of~$G$ (see~\cite[Corollary~19.13]{S}).

\medskip

Recall \cite{Dixmier} that a maximal abelian subalgebra $A$, or
masa, of a von Neumann algebra $M$ is called regular if its
normalizer $N(A)$ consisting of unitaries $u\in M$ such that
$uAu^*=A$ generates~$M$ as a von Neumann algebra, and singular if
the normalizer consists only of unitaries in $A$. If $A$ is
regular and there exists a faithful normal conditional expectation
of $M$ onto $A$ then $A$ is called Cartan~\cite{FM2}.

Since the action $T$ is free and ergodic, the algebras $C_T$ and
$S_T$ are maximal abelian in $R_T$. The algebra $C_T$ is Cartan.
Nielsen~\cite{Nielsen} was the first who noticed that if the
action is weakly mixing (i.e. the only eigenfunctions are
constants) then $S_T$ is singular (see~\cite{Packer2,SS} for
different proofs). More generally, the normalizer $N(S_T)$ always
depends only on the discrete part of the spectrum~\cite{Packer1}
(see also~\cite{Hahn}). We shall first give a short proof of this
result.

\begin{theorem} \label{2.1}
Let $L^\infty_0(X)$ be the subalgebra of $L^\infty(X)$ generated
by the eigenfunctions of the action~$\alpha$. Then the von Neumann
algebra $\NN(S_T)$ generated by $N(S_T)$ is
$L^\infty_0(X)\rtimes_\alpha G$.
\end{theorem}

\noindent{\it Proof.} If $u\in C_T$ is an eigenfunction,
$\alpha_g(u)=\<\chi,g\>u$ for some $\chi\in\hat G$, then since the
action is ergodic, $u$ is a unitary. It is in the normalizer of
$S_T$, $uv_gu^*=\<\chi,-g\>v_g$. Thus $L^\infty_0(X)\rtimes_\alpha
G\subset\NN(S_T)$.

Conversely, let $u\in N(S_T)$. Then $\Ad u$ defines an
automorphism of $S_T$ which corresponds to a measurable
transformation $\sigma$ of $\hat G$. Consider $R_T$ in the Hilbert
space $L^2(X)\otimes L^2(\hat G)$ as above. Then the operator
$v=u(1\otimes u_\sigma^*)$ commutes with $1\otimes L^\infty(\hat
G)$, hence it belongs to
$$
(L^\infty(X)\otimes B(L^2(\hat G)))\cap(1\otimes L^\infty(\hat
G))'=L^\infty(X)\otimes L^\infty(\hat G).
$$
Thus $v$ is given by a measurable family $\{v_\ell\}_{\ell\in\hat
G}$ of unitaries in $L^\infty(X)$. Since $u\in R_T$ and $v$
commutes with $1\otimes m_g$, we have
$$
u=(\alpha_g\otimes\Ad m_g^*)(u)=(\alpha_g\otimes\Ad
m_g^*)(v)(\alpha_g\otimes\Ad m_g^*)(1\otimes
u_\sigma)=(\alpha_g\otimes1)(v)(1\otimes m_g^*u_\sigma m_g).
$$
Hence $v=(\alpha_g\otimes1)(v)(1\otimes m_g^*u_\sigma
m_gu_\sigma^*)$. The operator $u_\sigma m_gu_\sigma^*$ is the
operator of multiplication by the function $g\circ\sigma^{-1}$.
Thus for almost all $\ell\in\hat G$
$$
v_\ell=\<\b{\ell}\sigma^{-1}(\ell),g\>\alpha_g(v_\ell).
$$
We see that for almost all $\ell$ the unitary $v_\ell$ lies in
$L^\infty_0(X)$, which means that $v\in L^\infty_0(X)\otimes
L^\infty(\hat G)$. Thus
$$
u=v(1\otimes u_\sigma)\in (L^\infty_0(X)\otimes B(L^2(\hat
G)))\cap R_T=(L^\infty_0(X)\otimes B(L^2(\hat
G)))^{\alpha\otimes\Ad m^*}=L^\infty_0(X)\rtimes_\alpha G.
$$
\endveriff

\noindent {\bf Remark.} The proof works without any
modifications in the case when a locally compact separable abelian
group acts ergodically on a von Neumann algebra with separable
predual.

\medskip

All the Cartan algebras $C_T$ are conjugate by a well-known result
of Dye~\cite{Dye2}, so the position of~$C_T$ inside~$R_T$ does not
contain any information about the original action. On the other
hand, the relative position of $C_T$ and $S_T$ defines the action.
More precisely, we have

\begin{proposition} \label{1.2}
Let $g\mapsto T^{(i)}_g\in\Aut(X_i,\mu_i)$ be a free
measure-preserving action of a countable abelian group~$G_i$,
$i=1,2$. Suppose there exists an isomorphism $\gamma\colon
R_{T^{(1)}}\to R_{T^{(2)}}$ such that
$\gamma(S_{T^{(1)}})=S_{T^{(2)}}$ and
$\gamma(C_{T^{(1)}})=C_{T^{(2)}}$. Then there exist an isomorphism
$S\colon(X_1,\mu_1)\to(X_2,\mu_2)$ of measure spaces and a group
isomorphism $\beta\colon G_2\to G_1$ such that
$T^{(2)}_g=ST^{(1)}_{\beta(g)}S^{-1}$ for $g\in G_2$.
\end{proposition}

\noindent{\it Proof.} The result follows easily from the fact that
the only unitaries in $S_T$ which normalize~$C_T$ are the scalar
multiples of $v_g$, $g\in G$. Indeed, if $v\in S_T$ normalizes
$C_T$ and $v=\sum_ga_gv_g$, $a_g\in\C$, is its Fourier series then
for arbitrary $x\in C_T$ the equality $vx=\alpha(x)v$ for $x\in
C_T$, where $\alpha=\Ad v$, implies $a_g\alpha_g(x)=a_g\alpha(x)$
for all $g\in G$. Thus $\alpha_g=\alpha$ if $a_g\ne0$. Since the
action is free, this means that $a_g\ne0$ for a unique $g$, and
$v=a_gv_g$. Hence if we have an isomorphism~$\gamma$ as in the
formulation of the proposition, then there exist an isomorphism
$\beta$ of $G_2$ onto $G_1$ and a character $\chi\in\hat G_2$ such
that $\gamma(v_{\beta(g)})=\<\chi,g\>v_g$ for $g\in G_2$. Then for
$x\in C_{T^{(1)}}$ and $g\in G_2$ we have
$\gamma(\alpha_{\beta(g)}(x))=\gamma(v_{\beta(g)}xv_{\beta(g)}^*)
=v_g\gamma(x)v_g^*=\alpha_g(\gamma(x))$. So for $S$ we can take
the transformation implementing the isomorphism $\gamma$ of
$C_{T^{(1)}}$ onto $C_{T^{(2)}}$.
\endverif

This observation leads to the following question. How much
information about the system is contained in the algebra $S_T$? If
the spectrum is purely discrete then $S_T$ is a Cartan subalgebra,
so in this case we get no information.

\medskip\noindent
{\bf Conjecture.} For weakly mixing systems the algebra $S_T$
determines the system completely. In other words, the assumption
$\gamma(C_{T^{(1)}})=C_{T^{(2)}}$ in Proposition~\ref{1.2} is
redundant.

\bigskip\bigskip

\section{Spectral invariants} \label{3}

One approach to the problem of conjugacy of masas in a
II$_1$-factor, initiated in the work of Pukanszky~\cite{Pukanszky}
is to consider together with a masa $A\subset M$ its conjugate
$JAJ$, where $J$ is the modular involution associated with a
tracial vector $\xi$, and then to consider the conjugacy problem
for such pairs in $B(L^2(M))$. We thus identify $A$ with an
algebra $L^\infty(Y,\nu)$ and consider a direct integral
decomposition of the representation $a\otimes b\mapsto aJb^*J$ of
the C$^*$-tensor product algebra $A\otimes A$. Thus we obtain a
measure class $[\eta]$ on $Y\times Y$ and a measurable field of
Hilbert spaces $\{H_{x,y}\}_{(x,y)\in Y\times Y}$ such that
$[\eta]$ is invariant with respect to the flip
$(x,y)\mapsto(y,x)$, its left (and right) projection onto $Y$ is
$[\nu]$, and
$$
L^2(M)=\int^\oplus_{Y\times Y}H_{x,y}d\eta(x,y),
$$
see~\cite{FM2} for details. Let $m(x,y)=\dim H_{x,y}$ be the
multiplicity function. Note that $m(x,x)=1$ and the subspace
$\int^\oplus_{Y\times Y}H_{x,x}d\eta(x,x)$ is identified with
$\overline{A\xi}$. Indeed, $\zeta\in L^2(M)$ lives on the diagonal
$\Delta(Y)\subset Y\times Y$ if and only if $a\zeta=Ja^*J\zeta$
for all $a\in A$. Since $A$ is maximal abelian, this is equivalent
to $\zeta\in\overline{A\xi}$. In particular, the projection
$e_A=[A\xi]$ corresponds to the characteristic function of
$\Delta(Y)$, so it belongs to $A\vee JAJ$ (see~\cite{Popa_notes}).

The triple $(Y,[\eta],m)$ is a conjugacy invariant for the pair
$(A,J)$ in the following sense. If $A\subset M$ and $B\subset N$
are masas then a unitary $U\colon L^2(M)\to L^2(N)$ such that
$UAU^*=B$ and $UJ_MU^*=J_N$ exists if and only if there exists an
isomorphism $F\colon(Y_A,[\nu_A])\to(Y_B,[\nu_B])$ such that
$(F\times F)_*([\eta_A])=[\eta_B]$ and $m_B\circ(F\times F)=m_A$.
Indeed, the fact that $U$ defines $F$ follows by definition.
Conversely, for given $F$ we can suppose without loss of
generality that $\eta_A$ is invariant with respect to the flip and
$(F\times F)_*(\eta_A)=\eta_B$. Then there exists a measurable
field of unitaries $\tilde U_{x,y}\colon H^A_{x,y}\to
H^B_{F(x),F(y)}$, and we can define the unitary $\tilde
U=\int^\oplus_{Y_A\times Y_A}\tilde U_{x,y}d\eta_A(x,y)$. It has
the property $\tilde UA\tilde U^*=B$. We want to modify $\tilde U$
in a way such that the condition $UJ_MU^*=J_N$ is also satisfied.
Note that $J_M$ is given by a measurable field of anti-unitaries
$J^A_{x,y}\colon H^A_{x,y}\to H^A_{y,x}$ such that
$J^A_{y,x}J^A_{x,y}=1$, and analogously $J_N$ defines a measurable
field $\{J^B_{x,y}\}_{x,y}$. We can easily arrange $\tilde
U_{x,x}J^A_{x,x}=J^B_{F(x),F(x)}\tilde U_{x,x}$. Outside of the
diagonal we choose a measurable subset $Z\subset Y_A\times Y_A$
which meets every two-point set $\{(x,y),(y,x)\}$ only once. Then
we define
$$
U_{x,y}=\cases{\tilde U_{x,y}\hspace{22mm}\ \ \hbox{if}\
(x,y)\in\Delta(Y_A)\cup Z,\cr J^B_{F(y),F(x)}\tilde
U_{y,x}J^A_{x,y}\ \ \hbox{otherwise}.}
$$
Then $U_{y,x}J^A_{x,y}=J^B_{F(x),F(y)}U_{x,y}$, so for
$U=\int^\oplus_{Y_A\times Y_A}U_{x,y}d\eta_A(x,y)$ we have
$UJ_M=J_NU$.

\smallskip

A rougher invariant is the set $P(A)\subset\N\cup\{\infty\}$ of
essential values of the multiplicity function~$m$ on $(Y\times
Y)\backslash\Delta(Y)$, which was introduced by
Pukanszky~\cite{Pukanszky} (we rather use the definition of
Popa~\cite{Popa_notes}). In other words, $P(A)$ is the set of $n$
such that the type I algebra $(A\vee JAJ)'(1-e_A)$ has a non-zero
component of type I$_n$. This invariant solves a weaker conjugacy
problem: $P(A)=P(B)$ if and only if there exists a unitary $U$
such that $U(A\vee J_MAJ_M)U^*=B\vee J_NBJ_N$ and $Ue_AU^*=e_B$.

\medskip

Return to our masas $S_T$ in $R_T$. As above, consider $R_T$
acting on $L^2(X\times\hat G)$ with the modular involution given
by~(\ref{e2.1}) and~(\ref{e3.1}). For the construction of the
triple $(Y_T,[\eta_T],m_T)$ for the masa~$S_T$ it is natural to
take $Y_T=\hat G$. Let $\mu_T$ and $n_T$ be the spectral measure
and the multiplicity function of the representation $g\mapsto
u_g$, so that
$$
L^2(X)=\int^\oplus_{\hat G}H_\ell\,d\mu_T(\ell).
$$
Following \cite{Hahn} we have a direct integral decomposition
$$
L^2(X\times\hat G)=\int^\oplus_{\hat G\times\hat G
}H_{\ell_2}\,d\lambda(\ell_1)d\mu_T(\ell_2),
$$
with respect to which $v_g=1\otimes m_g$ corresponds to the
function $(\ell_1,\ell_2)\mapsto g(\ell_1)$, while
$Jv_g^*J=u_{-g}\otimes m_g$ corresponds to $(\ell_1,\ell_2)\mapsto
g(\ell_1\b{\ell_2})$. Hence if we define $\eta_T$ as the image of
the measure $\lambda\times\mu_T$ under the map $\hat G\times\hat
G\to\hat G\times\hat G$,
$(\ell_1,\ell_2)\mapsto(\ell_1,\ell_1\b{\ell_2})$, then with
respect to the decomposition
$$
L^2(X\times\hat G)=\int^\oplus_{\hat G\times\hat G
}H_{\ell_1\b{\ell_2}}\,d\eta_T(\ell_1,\ell_2)
$$
the operator $v_g$ corresponds to the function
$(\ell_1,\ell_2)\mapsto g(\ell_1)$, while $Jv_g^*J$ corresponds to
$(\ell_1,\ell_2)\mapsto g(\ell_2)$. This is the decomposition we
are looking for. Thus we have proved the following (see
also~\cite{Hahn}).

\begin{proposition} \label{3.1}
The triple $(Y_T,[\eta_T],m_T)$ associated with the masa $S_T$ in
$R_T$ is given by $Y_T=\hat G$, $\int fd\eta_T=\int
f(\ell_1,\ell_1\b{\ell_2})d\lambda(\ell_1)d\mu_T(\ell_2)$,
$m_T(\ell_1,\ell_2)=n_T(\ell_1\b{\ell_2})$, where $\mu_T$ and
$n_T$ are the spectral measure and the multiplicity function for
the representation $g\mapsto u_g$ of $G$.
\endveriff
\end{proposition}

\begin{corollary}
The Pukanszky invariant $P(S_T)$ is the set of essential values of
the multiplicity function $n_T$ on $\hat G\backslash\{e\}$.
\endveriff
\end{corollary}

This corollary is also obvious from
$$
S_T\vee JS_TJ=\{u_g\ |\ g\in G\}''\otimes L^\infty(\hat G),\ \
e_{S_T}=p_1\otimes1,
$$
where $p_1\in B(L^2(X))$ is the projection onto the constants.

\smallskip

Pukanszky introduced his invariant to construct a countable family
of non-conjugate singular masas in the hyperfinite II$_1$-factor.
For each $n\in\N$ he constructed a singular masa $A$ with
$P(A)=\{n\}$. Thanks to advances in the spectral theory of
dynamical systems~\cite{KL} we now know much more.

\begin{corollary}
For any subset $E$ of $\N$ containing $1$ there exists a weakly
mixing automorphism~$T$ such that $P(S_T)=E$.
\endveriff
\end{corollary}

If the spectrum of the representation $g\mapsto u_g$ is Lebesgue,
i.e. the spectral measure $\mu_T$ is equivalent to the Haar
measure $\lambda$ on $\hat G\backslash\{e\}$, then
$[\eta_T]=[\lambda\times\lambda]$ on $(\hat G\times \hat
G)\backslash\Delta(\hat G)$. Hence if we have two such systems
then any measurable isomorphism $F\colon (\hat
G_1,[\lambda_1])\to(\hat G_2,[\lambda_2])$ has the property
$(F\times F)_*([\eta_{T^{(1)}}])=[\eta_{T^{(2)}}]$. Thus we have

\begin{corollary}
Let $g\mapsto T^{(i)}_g\in\Aut(X_i,\mu_i)$ be a free ergodic
measure-preserving action of a countable abelian group~$G_i$,
$i=1,2$. Suppose these actions have homogeneous Lebesgue spectra
of the same multiplicity. Then for any $*$-isomorphism
$\gamma\colon S_{T^{(1)}}\to S_{T^{(2)}}$ there exists a unitary
$U\colon L^2(R_{T^{(1)}})\to L^2(R_{T^{(2)}})$ such that
$UaU^*=\gamma(a)$ for $a\in S_{T^{(1)}}$ and
$UJ_{T^{(1)}}U^*=J_{T^{(2)}}$.
\endveriff
\end{corollary}

It is clear, however, that in order to be extended to an
isomorphism of $R_{T^{(1)}}$ on $R_{T^{(2)}}$, $\gamma$~has to be
at least trace-preserving. But even this is not always enough, see
Section~\ref{5}. Thus for such system as Bernoulli shifts, which
have countably multiple Lebesgue spectra, the invariant
$(Y_T,[\eta_T],m_T)$ does not contain any useful information.

\bigskip\bigskip

\section{The isomorphism problem} \label{4}

As a partial result towards a proof of our conjecture we have

\begin{theorem} \label{4.1}
Let $g\mapsto T^{(i)}_g\in\Aut(X_i,\mu_i)$ be a weakly mixing free
measure-preserving action of a countable abelian group~$G_i$,
$i=1,2$. Suppose there exists an isomorphism $\gamma\colon
R_{T^{(1)}}\to R_{T^{(2)}}$ such that
$\gamma(S_{T^{(1)}})=S_{T^{(2)}}$ and such that the Cartan
algebras $\gamma(C_{T^{(1)}})$ and $C_{T^{(2)}}$ are inner
conjugate in~$R_{T^{(2)}}$. Then there exist an isomorphism
$S\colon(X_1,\mu_1)\to(X_2,\mu_2)$ of measure spaces and a group
isomorphism $\beta\colon G_2\to G_1$ such that
$T^{(2)}_g=ST^{(1)}_{\beta(g)}S^{-1}$ for $g\in G_2$.
\end{theorem}

We shall also describe explicitly all possible isomorphisms
$\gamma$ as in the theorem. In other words, for a weakly mixing
free measure-preserving action $T$ of a countable abelian group
$G$ on $(X,\mu)$ we shall compute the group $\AutT$ consisting of
all automorphisms $\gamma$ of $R_T$ with the properties
$\gamma(S_T)=S_T$, and the masas $\gamma(C_T)$ and $C_T$ are inner
conjugate.

Recall (see \cite{FM2}) that any automorphism $S$ of the orbit
equivalence relation defined by the action of $G$ extends
canonically to an automorphism $\alpha_S$ of $R_T$. Such an
automorphism leaves~$S_T$ invariant if and only if there exists an
automorphism $\beta$ of $G$ such that $T_gS=ST_{\beta(g)}$. Denote
by $I(T)$ the group of all such transformations $S$. For $S\in
I(T)$, $\alpha_S$ is defined by the equalities $\alpha_S(f)=f\circ
S^{-1}$  for $f\in C_T=L^\infty(X,\mu)$,
$\alpha_S(v_g)=v_{\beta^{-1}(g)}$ for $g\in G$. Consider also the
dual action $\sigma$ of $\hat G$ on $R_T$, $\sigma_\chi(f)=f$ for
$f\in C_T$, $\sigma_\chi(v_g)=\<\chi,-g\>v_g$. The group of
automorphisms of the form $\sigma_\chi\circ\alpha_S$ ($\chi\in\hat
G$ and $S\in I(T)$) is the intersection of the groups
$\Aut(R_T,C_T)$ and $\Aut(R_T,S_T)$. It turns out that up to inner
automorphisms defined by unitaries in $S_T$ such automorphisms
exhaust the whole group $\AutT$.

\begin{theorem} \label{4.2}
The group $\AutT$ of automorphisms $\gamma$ of $R_T$ for which
$\gamma(S_T)=S_T$, and $\gamma(C_T)$ and $C_T$ are inner
conjugate, consists of elements of the form $\Ad
w\circ\sigma_\chi\circ\alpha_S$, where $w\in S_T$, $\chi\in\hat
G$, $S\in I(T)$.
\end{theorem}

We conjecture that in fact this theorem gives the description of
the group~$\Aut(R_T,S_T)$.

\smallskip

It is well-known that all Cartan subalgebras of the hyperfinite
II$_1$-factor are conjugate~\cite{CFS}, so they are approximately
inner conjugate in an appropriate sense. It is known also that if
the $L^2$-distance between the unit balls of two Cartan
subalgebras is less than one, then they are inner
conjugate~\cite{Popa_countable,Popa_unit}. However, there exists
an uncountable family of Cartan subalgebras, no two of which are
inner conjugate~\cite{Packer1}.

\smallskip

We shall first prove that Theorem~\ref{4.1} follows from
Theorem~\ref{4.2}.
Consider the group $G=G_1\times G_2$ and its action $T$ on
$(X,\mu)=(X_1\times X_2,\mu_1\times\mu_2)$,
$T_{(g_1,g_2)}=T^{(1)}_{g_1}\times T^{(2)}_{g_2}$. Then $R_T$ can
be identified with $R_{T^{(1)}}\otimes R_{T^{(2)}}$ in such a way
that $C_T=C_{T^{(1)}}\otimes C_{T^{(2)}}$,
$v_{(g_1,g_2)}=v_{g_1}\otimes v_{g_2}$. Consider the automorphism
$\tilde\gamma$ of $R_T$,
$$
\tilde\gamma(a\otimes b)=\gamma^{-1}(b)\otimes\gamma(a).
$$
By Theorem \ref{4.2}, $\tilde\gamma$ must be of the form $\Ad
w\circ\sigma_\chi\circ\alpha_{\tilde S}$ with $w\in S_T$,
$\chi=(\chi_1,\chi_2)\in\hat G_1\times\hat G_2$ and $\tilde S\in
I(T)$. Let $\tilde\beta\in\Aut(G)$ be such that $T_g\tilde
S=\tilde ST_{\tilde\beta(g)}$. Since $\tilde\gamma^2=\id$, we have
$\tilde\beta^2=\id$. Define the homomorphism $\beta\colon G_2\to
G_1$ as the composition of the map $g_2\mapsto\tilde\beta(0,g_2)$
with the projection $G_1\times G_2\to G_1$, and $\beta'\colon
G_1\to G_2$ as the composition of the map
$g_1\mapsto\tilde\beta(g_1,0)$ with the projection $G_1\times
G_2\to G_2$. Fix $g_2\in G_2$. Then
$\tilde\beta(0,g_2)=(\beta(g_2),h)$ for some $h\in G_2$. We have
$$
\gamma^{-1}(v_{g_2})\otimes1=\tilde\gamma(1\otimes v_{g_2})
=\<\chi_1,-\beta(g_2)\>\<\chi_2,-h\>v_{\beta(g_2)}\otimes v_h.
$$
It follows that $h=0$, that is
$\tilde\beta(0,g_2)=(\beta(g_2),0)$. Analogously
$\tilde\beta(g_1,0)=(0,\beta'(g_1))$. Thus
$\tilde\beta(g_1,g_2)=(\beta(g_2),\beta'(g_1))$. Since
$\tilde\beta^2=\id$, we conclude that $\beta'=\beta^{-1}$. Then the
identity $T_g\tilde S=\tilde ST_{\tilde\beta(g)}$ is rewritten in
terms of the actions on $L^\infty(X_1\times X_2)$ as
$$
(\alpha_{g_1}\otimes\alpha_{g_2})\circ\alpha_{\tilde
S}=\alpha_{\tilde
S}\circ(\alpha_{\beta(g_2)}\otimes\alpha_{\beta^{-1}(g_1)}).
$$
Letting $g_2=0$ we see that for $f\in L^\infty(X_1)$
$$
((\alpha_{g_1}\otimes1)\circ\alpha_{\tilde
S})(f\otimes1)=(\alpha_{\tilde
S}\circ(1\otimes\alpha_{\beta^{-1}(g_1)}))(f\otimes1)=\alpha_{\tilde
S}(f\otimes1),
$$
so that $\alpha_{\tilde S}(L^\infty(X_1)\otimes1)\subset
L^\infty(X_1\times X_2)^{\alpha_{G_1}\otimes1}=1\otimes
L^\infty(X_2)$. Analogously $\alpha_{\tilde S}(1\otimes
L^\infty(X_2))\subset L^\infty(X_1)\otimes1$. It follows that
$$
\alpha_{\tilde S}(L^\infty(X_1)\otimes1)=1\otimes L^\infty(X_2)\ \
\hbox{and}\ \ \alpha_{\tilde S}(1\otimes
L^\infty(X_2))=L^\infty(X_1)\otimes1.
$$
Hence there exist isomorphisms $S\colon(X_1,\mu_1)\to(X_2,\mu_2)$
and $S'\colon(X_2,\mu_2)\to(X_1,\mu_1)$ such that for almost all
$(x_1,x_2)$ we have $\tilde S(x_1,x_2)=(S'x_2,Sx_1)$. The identity
$(T^{(1)}_{g_1}\times T^{(2)}_{g_2})\tilde S=\tilde S
(T^{(1)}_{\beta(g_2)}\times T^{(2)}_{\beta^{-1}(g_1)})$ implies
that $T^{(2)}_{g_2}S=ST^{(1)}_{\beta(g_2)}$.

\medskip

Now we turn to the proof of Theorem \ref{4.2}. The proof will be
given in a series of lemmas. Let $\gamma\in\AutT$.

\begin{lemma}
The automorphism $\gamma$ can be implemented by a unitary $U$ on
$L^2(R_T)$ such that $UJC_TJU^*=JC_TJ$, where $J$ is the modular
involution.
\end{lemma}

\noindent{\it Proof.} Let $\tilde U$ be the canonical
implementation of $\gamma$ commuting with $J$. From the assumption
that $C_T$ and $\gamma(C_T)$ are inner conjugate we can choose
$u\in R_T$ such that $uC_Tu^*=\gamma(C_T)$. Then we can take
$U=Ju^*J\tilde U$.
\endverif

Representing $R_T$ on $L^2(X)\otimes L^2(\hat G)$ as usual, so
that $JC_TJ=L^\infty(X)\otimes1$ and $S_T=1\otimes L^\infty(\hat
G)$ (see~(\ref{e3.1})), we conclude that $\Ad U$ defines
measure-preserving transformations $S_1$ of $X$ and $\sigma$ of
$\hat G$. Then $W=U(u_{S_1}^*\otimes u_\sigma^*)$ commutes with
$L^\infty(X)\otimes1$ and $1\otimes L^\infty(\hat G)$, hence it is
a unitary in $L^\infty(X\times\hat G)$. For $\ell\in\hat G$ denote
by $w_\ell$ the function in $L^\infty(X)$ defined by
$w_\ell(x)=W(x,\ell)$.

Since $U$ defines an automorphism of $R_T$, for $f\in L^\infty(X)$
the element $U\pi(f)U^*$ must by~(\ref{e3.1}) commute with
$u_h\otimes m_h^*$.

\begin{lemma} \label{4.4}
With the above notations, for $\zeta\in L^2(\hat G,L^2(X))\cong
L^2(X)\otimes L^2(\hat G)$ we have
$$
(U\pi(f)U^*\zeta)(\ell)=\sum_{g\in
G}(\alpha_{S_1}\circ\alpha_g)(f)
 \int_{\hat G}\<\b{\sigma^{-1}(\ell)}\sigma^{-1}(\ell_1),g\>w_\ell
   w_{\ell_1}^*\zeta(\ell_1)d\lambda(\ell_1),
$$
$\displaystyle ((u_h\otimes m_h^*)U\pi(f)U^*(u_h^*\otimes
m_h)\zeta)(\ell)$
$$
 =\sum_{g\in G}(\alpha_h\circ\alpha_{S_1}\circ\alpha_g)(f)\int_{\hat
G}\<\b{\ell}\ell_1,h\>
\<\b{\sigma^{-1}(\ell)}\sigma^{-1}(\ell_1),g\> \alpha_h(w_\ell
w_{\ell_1}^*)\zeta(\ell_1)d\lambda(\ell_1).
$$
The above series are meaningless for fixed $\ell$ and should be
considered as series of functions in $L^2(X\times\hat G)$.
\end{lemma}

\noindent{\it Proof.} Note that
$(W\zeta)(\ell)=w_\ell\zeta(\ell)$, $((1\otimes
m_g)\zeta)(\ell)=\<\ell,g\>\zeta(\ell)$. The operator $u_\sigma
e_g u_\sigma^*$ is the projection onto the one-dimensional space spanned by the 
function $u_\sigma g\in
L^2(\hat G)$, so for $f\in L^2(\hat G)$,
$$
(u_\sigma e_g u_\sigma^*f)(\ell)=(u_\sigma g)(\ell)\cdot(f,u_\sigma g)
 =\int_{\hat G}\<\sigma^{-1}(\ell)
     \b{\sigma^{-1}(\ell_1)},g\>f(\ell_1)d\lambda(\ell_1).
$$
Hence $\displaystyle ((1\otimes u_\sigma e_g u_\sigma^*)\zeta)(\ell)
=\int_{\hat G}\<\sigma^{-1}(\ell)
     \b{\sigma^{-1}(\ell_1)},g\>\zeta(\ell_1)d\lambda(\ell_1)$.
Now we compute:
\begin{eqnarray*}
(U\pi(f)U^*\zeta)(\ell)
 &=&\left(W(u_{S_1}\otimes u_\sigma)
      \left(\sum_g\alpha_g(f)\otimes e_{-g}\right)
         (u_{S_1}^*\otimes u_\sigma^*)W^*\zeta\right)(\ell)\\
 &=&\sum_g(W((\alpha_{S_1}\circ\alpha_g)(f)
      \otimes u_\sigma e_{-g}u_\sigma^*)W^*\zeta)(\ell)\\
 &=&\sum_g w_\ell(\alpha_{S_1}\circ\alpha_g)(f)
      \int_{\hat G}\<\sigma^{-1}(\ell)\b{\sigma^{-1}(\ell_1)},-g\>
        (W^*\zeta)(\ell_1)d\lambda(\ell_1)\\
 &=&\sum_g (\alpha_{S_1}\circ\alpha_g)(f)
      \int_{\hat G}\<\b{\sigma^{-1}(\ell)}\sigma^{-1}(\ell_1),g\>
        w_\ell w_{\ell_1}^*\zeta(\ell_1)d\lambda(\ell_1),
\end{eqnarray*}
and \newline $\displaystyle ((u_h\otimes
m_h^*)U\pi(f)U^*(u_h^*\otimes m_h)\zeta)(\ell)$
\begin{eqnarray*}
 &=&\<\b{\ell},h\>u_h(U\pi(f)U^*(u_h^*\otimes m_h)\zeta)(\ell)\\
 &=&\<\b{\ell},h\>u_h\sum_g (\alpha_{S_1}\circ\alpha_g)(f)
      \int_{\hat G}\<\b{\sigma^{-1}(\ell)}\sigma^{-1}(\ell_1),g\>
        w_\ell w_{\ell_1}^*((u_h^*\otimes m_h)\zeta)(\ell_1)
          d\lambda(\ell_1)\\
 &=&\<\b{\ell},h\>u_h\sum_g (\alpha_{S_1}\circ\alpha_g)(f)
      \int_{\hat G}\<\b{\sigma^{-1}(\ell)}\sigma^{-1}(\ell_1),g\>
        w_\ell w_{\ell_1}^*u_h^*\<\ell_1,h\>\zeta(\ell_1)
          d\lambda(\ell_1)\\
 &=&\sum_g (\alpha_h\circ\alpha_{S_1}\circ\alpha_g)(f)
      \int_{\hat G}\<\b{\ell}\ell_1,h\>
        \<\b{\sigma^{-1}(\ell)}\sigma^{-1}(\ell_1),g\>
          \alpha_h(w_\ell w_{\ell_1}^*)\zeta(\ell_1)d\lambda(\ell_1).
\end{eqnarray*}
\endveriff

\begin{lemma} \label{4.5}
Let $g\mapsto P_g\in\AutX$ be a free measure-preserving action of
$G$, $Q\in\AutX$, $H$ a Hilbert space, $a_g$ and $b_g$ maps from
$X$ to $H$ such that
\newline
{\rm(i)} the vectors $a_g(x)$, $g\in G$, are mutually orthogonal
for almost all $x\in X$;
\newline
{\rm(ii)} $\sum_g||a_g(x)||^2$ is finite and non-zero for almost
all $x$;
\newline
and the same conditions hold for $\{b_g\}_g$. Suppose for all
$f\in L^\infty(X)$ and almost all $x\in X$
$$
\sum_g(\alpha_Q\circ\alpha_{P_g})(f)(x)a_g(x)
 =\sum_g\alpha_{P_g}(f)(x)b_g(x).
$$
Then $Q$ is in the full group generated by $P_g$, $g\in G$, and if
$g(x)\in G$ is such that $Q^{-1}x=P_{-g(x)}x$ then
$a_g(x)=b_{g+g(x)}(x)$ for all $g\in G$ and almost all $x\in X$.
\end{lemma}

\noindent{\it Proof.} Let $X_0=\{x\in X\,|\, Q^{-1}x\notin P_Gx,\,
P_gx\ne x\ \hbox{for}\ g\ne0\}$. There exists a countable family
$\{X_i\}_{i\in I}$ of measurable subsets of $X$ such that for
arbitrary finite subset $F$ of $G$ and almost all $x\in X_0$ there
exists $i\in I$ such that $x\in X_i$, the sets $P_gX_i$, $g\in F$,
are mutually disjoint and $Q^{-1}x\notin\cup_{g\in F}P_gX_i$.
Indeed, first note that choosing an arbitrary $Q$- and
$P_g$-invariant norm-separable weakly dense C$^*$-subalgebra $A$
of $L^\infty(X)$, we can identify the measure space $(X,\mu)$ with
the spectrum of~$A$. Thus without loss of generality we can
suppose that $X$ is a compact metric space and $Q$ and $P_g$ are
homeomorphisms. Moreover, by regularity of the measure it is
enough to prove the assertion for arbitrary compact subset $K$ of
$X_0$. But then for fixed~$F$ we can consider for each $x\in K$ a
neighborhood $U_x$ such that $P_gU_x$, $g\in F$, are disjoint,
$Q^{-1}U_x\cap P_gU_x=\varnothing$ for $g\in F$, and then choose a
finite subcovering from $\{U_x\}_{x\in K}$.

Consider the countable set ${\cal F}\subset L^\infty(X)$
consisting of characteristic functions of the sets $X_i$, $i\in
I$, and all their translations under the action of $G$. For almost
all $x\in X_0$ and all $f\in{\cal F}$ the assumptions of the lemma
are satisfied. Let $x\in X_0$ be such a point. Fix $h\in G$. For
arbitrary finite subset $F$ of $G$, $h\in F$, there exists
$f\in{\cal F}$ such that $\alpha_{P_h}(f)(x)=1$,
$\alpha_{P_g}(f)(x)=0$ for $g\in F\backslash\{h\}$ and
$(\alpha_Q\circ\alpha_{P_g})(f)(x)=0$ for $g\in F$. Then
\begin{eqnarray*}
||b_h(x)||
 & = &\left\|\sum_{g\notin F}(\alpha_Q\circ\alpha_{P_g})(f)(x)a_g(x)
        -\sum_{g\notin F}\alpha_{P_g}(f)(x)b_g(x)\right\|\\
 &\le&\left(\sum_{g\notin F}||a_g(x)||^2\right)^{1/2}
        +\left(\sum_{g\notin F}||b_g(x)||^2\right)^{1/2}.
\end{eqnarray*}
It follows that $b_h(x)=0$. But this contradicts the assumption
$\sum_h ||b_h(x)||^2>0$. Hence the set~$X_0$ has zero measure.
Thus $Q$ is indeed in the full group generated by $P_g$.

Let $Q^{-1}x=P_{-g(x)}x$. In the same way as above (or by
referring to the Rokhlin lemma) we can find a countable collection
${\cal F}$ of characteristic functions such that for almost all
$x\in X$ and arbitrary finite $F\subset G$, $0\in F$, there exists
$f\in{\cal F}$ such that $f(x)=1$, $\alpha_{P_g}(f)(x)=0$ for
$g\in F\backslash\{0\}$. Then
$$
a_{-g(x)}(x)-b_0(x)=\sum_{g\notin F}\alpha_{P_g}(f)(x)b_g(x)
 -\sum_{g\notin F-g(x)}(\alpha_Q\circ\alpha_{P_g})(f)(x)a_g(x),
$$
and we conclude that $a_{-g(x)}(x)=b_0(x)$. Replacing $f$ by
$\alpha_{P_h}(f)$ in the formulation of the lemma we see that its
assumptions are also satisfied for the collections $\{a_{g-h}\}_g$
and $\{b_{g-h}\}_g$, so that $a_{-g(x)-h}(x)=b_{-h}(x)$.
\endverif

Fix $h\in G$ and apply Lemma~\ref{4.5} to $P_g=S_1T_gS_1^{-1}$,
$Q=T_h$, $H=L^2(\hat G)$,
\begin{eqnarray*}
a_g(x)(\ell)
 &=&\int_{\hat G}\<\b{\ell}\ell_1,h\>
     \<\b{\sigma^{-1}(\ell)}\sigma^{-1}(\ell_1),g\>
       \alpha_h(w_\ell w_{\ell_1}^*)(x)d\lambda(\ell_1)\\
b_g(x)(\ell)
 &=&\int_{\hat G}
     \<\b{\sigma^{-1}(\ell)}\sigma^{-1}(\ell_1),g\>(w_\ell
       w_{\ell_1}^*)(x)d\lambda(\ell_1).
\end{eqnarray*}
To see that the assumptions of the lemma are satisfied, note that
up to the factor $\ell\mapsto\<\b{\ell},h\>\alpha_h(w_\ell)(x)$
the series $\sum_ga_g(x)$ is the Fourier series of the function
$\ell\mapsto\<\ell,h\>\alpha_h(w_\ell^*)(x)$ with respect to the
orthonormal basis $\{\b{u_\sigma g}\}_{g\in G}$.

Thus by Lemmas~\ref{4.4} and~\ref{4.5} we conclude that there
exists $g(h,x)$ such that $T_{-h}x=S_1T_{-g(h,x)}S_1^{-1}x$ and
$a_g(x)=b_{g+g(h,x)}(x)$, that is
$$
\int_{\hat G}\<\sigma^{-1}(\ell_1),g\>\Bigl(\<\b{\ell}\ell_1,h\>
\alpha_h(w_\ell w_{\ell_1}^*)(x)-
         \<\b{\sigma^{-1}(\ell)}\sigma^{-1}(\ell_1),g(h,x)\>(w_\ell
       w_{\ell_1}^*)(x)\Bigr)d\lambda(\ell_1)=0.
$$
Since the functions $u_\sigma
g=(\ell_1\mapsto\<\sigma^{-1}(\ell_1),g\>)$, $g\in G$, form an
orthonormal basis of $L^2(\hat G)$, we conclude that for almost
all $(x,\ell,\ell_1)$
\begin{equation} \label{e4.1}
\alpha_h(w_\ell w_{\ell_1}^*)(x)=\<\ell\b{\ell_1},h\>
\<\b{\sigma^{-1}(\ell)}\sigma^{-1}(\ell_1),g(h,x)\> (w_\ell
w_{\ell_1}^*)(x).
\end{equation}

\begin{lemma} \label{4.6}
There exists a continuous automorphism $\sigma_0$ of $\hat G$ and
$\chi\in\hat G$ such that $\sigma(\ell)=\chi\sigma_0(\ell)$ for
almost all $\ell$.
\end{lemma}

\noindent{\it Proof.} Replace $\ell$ by $\sigma(\ell)$ and
$\ell_1$ by $\sigma(\ell_1)$ in (\ref{e4.1}). Then we get
\begin{equation} \label{e4.2}
\alpha_h(w_{\sigma(\ell)}w_{\sigma(\ell_1)}^*)(x)
=\<\sigma(\ell)\b{\sigma(\ell_1)},h\>\<\b{\ell}\ell_1,g(h,x)\>
(w_{\sigma(\ell)}w_{\sigma(\ell_1)}^*)(x).
\end{equation}
Now substitute $\ell\ell_2$ for $\ell$ and $\ell_1\ell_2$ for
$\ell_1$. We get
\begin{equation} \label{e4.3}
\alpha_h(w_{\sigma(\ell\ell_2)}w_{\sigma(\ell_1\ell_2)}^*)(x)
=\<\sigma(\ell\ell_2)\b{\sigma(\ell_1\ell_2)},h\>
   \<\b{\ell}\ell_1,g(h,x)\>
   (w_{\sigma(\ell\ell_2)}w_{\sigma(\ell_1\ell_2)}^*)(x).
\end{equation}
Multiplying (\ref{e4.2}) by the equation conjugate to (\ref{e4.3})
we see that for almost all $(\ell,\ell_1,\ell_2)$ the element
$w_{\sigma(\ell)}w_{\sigma(\ell_1)}^*
w_{\sigma(\ell_1\ell_2)}w_{\sigma(\ell\ell_2)}^*$ is an
eigenfunction with eigenvalue $\sigma(\ell)\b{\sigma(\ell_1)}
\sigma(\ell_1\ell_2)\b{\sigma(\ell\ell_2)}$. Since the action is
weakly mixing, we conclude that
$$
\sigma(\ell\ell_2)\b{\sigma(\ell)}
=\sigma(\ell_1\ell_2)\b{\sigma(\ell_1)}
$$
(this is the only place where we use weak mixing instead of
ergodicity). Hence there exists a measurable map $\tilde\sigma_0$
of~$\hat G$ onto itself such that
$\tilde\sigma_0(\ell_2)=\sigma(\ell\ell_2)\b{\sigma(\ell)}$ for
almost all $(\ell,\ell_2)$. Then for almost all $(\ell_1,\ell_2)$
$$
\tilde\sigma_0(\ell_1\ell_2)
=\sigma(\ell\ell_1\ell_2)\b{\sigma(\ell)}
=\sigma(\ell\ell_1\ell_2)\b{\sigma(\ell\ell_2)}
  \sigma(\ell\ell_2)\b{\sigma(\ell)}
=\tilde\sigma_0(\ell_1)\tilde\sigma_0(\ell_2).
$$
So $\tilde\sigma_0$ is essentially a homomorphism, and since it is
measurable, it coincides almost everywhere with a continuous
homomorphism $\sigma_0$. Choose a character $\ell_1$ such that the
equality $\sigma_0(\ell)=\sigma(\ell_1\ell)\b{\sigma(\ell_1)}$
holds for almost all $\ell$. Set
$\chi=\sigma(\ell_1)\b{\sigma_0(\ell_1)}$. Then
$\sigma(\ell)=\chi\sigma_0(\ell)$ for almost all $\ell$. Since
$\sigma$ is an invertible measure-preserving transformation,
$\sigma_0$ must be an automorphism.
\endverif

Now we can rewrite (\ref{e4.1}) as
\begin{equation} \label{e4.4}
\alpha_h(w_\ell w_{\ell_1}^*)(x)
=\<\ell\b{\ell_1},h\>\<\sigma_0^{-1}(\b{\ell}\ell_1),g(h,x)\>
(w_\ell w_{\ell_1}^*)(x).
\end{equation}

\begin{lemma} \label{4.7}
Let $\ell_1$ be such that (\ref{e4.4}) holds for almost all
$(x,\ell)\in X\times\hat G$. Then there exist a unitary $b$ in
$L^\infty(\hat G)$ and a measurable map $e\colon X\to G$ such that
for almost all $(x,\ell)$ we have
$$
(w_{\ell\ell_1}w_{\ell_1}^*)(x)=b(\ell)\<\ell,e(x)\>.
$$
For all $h\in G$ and almost all $x\in X$ we have
$$
g(h,x)=\beta(e(x))-\beta(e(T_{-h}x))+\beta(h),
$$
where $\beta$ is the automorphism of $G$ dual to $\sigma_0$, i.e.
$\<\sigma_0(\ell),g\>=\<\ell,\beta(g)\>$.
\end{lemma}

\noindent{\it Proof.} Denote $w_{\ell\ell_1}w_{\ell_1}^*$ by
$v_\ell$. Then by (\ref{e4.4})
\begin{equation} \label{e4.5}
\alpha_h(v_\ell)(x) =\<\ell,h\>\<\sigma_0^{-1}(\b{\ell}),g(h,x)\>
(v_\ell)(x).
\end{equation}
Multiplying these identities for $v_\ell$, $v_{\ell_2}$ and
$v_{\ell\ell_2}^*$ we see that the function $c(\ell,\ell_2)=v_\ell
v_{\ell_2}v_{\ell\ell_2}^*$ is $G$-invariant, so it is a constant.
Thus we obtain a measurable symmetric (i.e.
$c(\ell,\ell_2)=c(\ell_2,\ell)$) 2-cocycle on $G$ with values in
$\T$. Since $G$ is abelian, any such a cocycle is a coboundary
(see e.g. \cite{Moore4}),
$c(\ell,\ell_2)=b(\ell)b(\ell_2)\b{b(\ell\ell_2)}$. Then
$\ell\mapsto\b{b(\ell)}v_\ell$ is a measurable homomorphism of
$\hat G$ into the unitary group of $L^\infty(X)$. By
\cite[Theorem~1]{Moore4} there exists a measurable map $e\colon
X\to G$ such that $\b{b(\ell)}v_\ell(x)=\<\ell,e(x)\>$.

Equation (\ref{e4.5}) implies that
$$
\<\ell,e(T_{-h}x)\>=\<\ell,h\>\<\sigma_0^{-1}(\b{\ell}),g(h,x)\>
\<\ell,e(x)\>,
$$
equivalently,
$$
\<\ell,e(T_{-h}x)-h+\beta^{-1}(g(h,x))-e(x)\>=1,
$$
from what the second assertion of the lemma follows.
\endverif

Recall that $S_1$ is the transformation of $X$ defined by $\Ad
U|_{L^\infty(X)}$.

\begin{lemma} \label{4.8}
Define a measurable map $S_2$ of $X$ onto itself by letting
$$
S_2x=S_1T_{-\beta(e(x))}S_1^{-1}x.
$$
Then $S_2$ is invertible and measure-preserving. Its inverse is
given by
$$
S_2^{-1}x=T_{e(x)}x.
$$
\end{lemma}

\noindent{\it Proof.} Recall that $g(h,x)$ was defined by the
equality $T_{-h}x=S_1T_{-g(h,x)}S_1^{-1}x$. Since by
Lemma~\ref{4.7}, $g(-e(x),x)=-\beta(e(T_{e(x)}x))$, it follows
that
$$
S_2T_{e(x)}x =S_1T_{-\beta(e(T_{e(x)}x))}S_1^{-1}T_{e(x)}x
=S_1T_{-\beta(e(T_{e(x)}x))-g(-e(x),x)}S_1^{-1}x=x.
$$
Hence $S_2$ is essentially surjective. Since it is also one-to-one
and measure-preserving on the sets $e^{-1}(\{g\})$, we conclude
that $S_2$ is invertible, measure-preserving and its inverse is
given by $S_2^{-1}x=T_{e(x)}x$.
\endverif

The final step is

\begin{lemma} \label{4.9}
The mapping $S=S_2^{-1}S_1$ has the property $T_gS=ST_{\beta(g)}$.
\end{lemma}

\noindent{\it Proof.} We compute:
\begin{eqnarray*}
S^{-1}T_{-h}x
 &=&S_1^{-1}S_2T_{-h}x\\
 &=&S_1^{-1}S_1T_{-\beta(e(T_{-h}x))}S_1^{-1}T_{-h}x\\
 &=&T_{-\beta(e(T_{-h}x))-g(h,x)}S_1^{-1}x\\
 &=&T_{-\beta(h)-\beta(e(x))}S_1^{-1}x\\
 &=&T_{-\beta(h)}S_1^{-1}S_2x\\
 &=&T_{-\beta(h)}S^{-1}x,
\end{eqnarray*}
where in the fourth equality we used Lemma~\ref{4.7}.
\endverif

Summarizing the results of Lemmas~\ref{4.6}-\ref{4.9}, we can
decompose $U=W(u_{S_1}\otimes u_\sigma)$ as follows. First, by
Lemma~\ref{4.7} for almost all $(x,\ell)$,
$w_\ell(x)=\<\ell\b{\ell_1},e(x)\>w_{\ell_1}(x)b(\ell\b{\ell_1})$,
so $W$ is the product of $u'\otimes1$, $v$ and $1\otimes w$, where
$u'\in L^\infty(X)$, $u'(x)=\<\b{\ell_1},e(x)\>w_{\ell_1}(x)$,
$v\in L^\infty(X\times\hat G)$, $v(x,\ell)=\<\ell,e(x)\>$, and
$w\in L^\infty(\hat G)$, $w(\ell)=b(\ell\b{\ell_1})$. By
Lemmas~\ref{4.8} and~\ref{4.9}, $u_{S_1}=u_{S_2}u_S$. Finally by
Lemma~\ref{4.6}, $u_\sigma=\lambda_\chi u_{\sigma_0}$, where
$\lambda_\chi$ is the operator of the left regular representation
of $\hat G$ on~$L^2(\hat G)$. Thus with $v'=v(u_{S_2}\otimes1)$ we
have
$$
U=(u'\otimes1)v(1\otimes
w)(u_{S_2}u_S\otimes1)(1\otimes\lambda_\chi u_{\sigma_0})
=(u'\otimes1)v'(1\otimes w)(1\otimes\lambda_{\chi})(u_S\otimes
u_{\sigma_0}).
$$
The unitaries $u'\otimes1$ and $v'$ both lie in the
commutant~$R_T'$. This is obvious for $u'\otimes1$ and follows for
$v'$ from the formula
$$
v'=\sum_g(p_g\otimes1)(u_g^*\otimes m_g),
$$
where $p_g$ is the characteristic function of the set
$e^{-1}(\{g\})$. Indeed, if $x\in e^{-1}(\{g\})$ then
$S_2^{-1}x=T_gx$ by Lemma~\ref{4.8}, and hence for arbitrary
$\zeta\in L^2(X\times\hat G)$ we have
$$
((p_g\otimes1)(u_g^*\otimes
m_g)\zeta)(x,\ell)=\<\ell,g\>\zeta(T_gx,\ell)
=\<\ell,e(x)\>\zeta(S_2^{-1}x,\ell) =(v'\zeta)(x,\ell).
$$
Thus the automorphism $\gamma$ is implemented by the unitary
$(1\otimes w)(1\otimes\lambda_{\chi})(u_S\otimes u_{\sigma_0})$,
so $\gamma=\Ad w\circ\sigma_\chi\circ\alpha_S$, and the proof of
Theorem~\ref{4.2} is complete.

\bigskip

From the definition of the group $\AutT$ it is unclear whether it
is a closed subgroup of $\Aut(R_T,S_T)$ (in the topology of
point-wise strong convergence). But if our conjecture that this
group coincides with $\Aut(R_T,S_T)$ (which is stronger than our
main conjecture in Section~\ref{2}) is true, then this group must
be closed. We shall prove that it is closed under slightly
stronger assumptions than weak mixing.

Recall that an action $T$ is called {\it rigid} if there exists a
sequence $\{g_n\}_n$ such that $g_n\to\infty$ and $u_{g_n}\to1$
strongly.

\begin{proposition}
Suppose $T$ is a weakly mixing action which is not rigid. Then the
group $\AutT$ is closed in $\Aut(R_T)$.
\end{proposition}

\noindent{\it Proof.} Suppose a sequence
$\{\alpha_n\}_n\subset\AutT$ converges to an automorphism
$\alpha$. By Theorem~\ref{4.2}, $\alpha_n=\sigma_{\chi_n}\circ\Ad
w_n\circ\alpha_{S_n}$. Passing to a subsequence we may suppose
that the sequence $\{\chi_n\}_n$ converges to a character $\chi$.
Then $\{\Ad w_n\circ\alpha_{S_n}\}_n$ converges to
$\sigma_\chi^{-1}\circ\alpha$, so to simplify the notations we may
suppose that all characters $\chi_n$ are trivial.

Let $f$ be a unitary generating $C_T$. Set
$$
\delta=\inf_{g\ne0}\|\alpha_g(f)-f\|_2.
$$
Since the action is not rigid, $\delta>0$. Suppose for some $n$
and $m$
$$
\|(\Ad w_n\circ\alpha_{S_n})(f)-(\Ad
w_m\circ\alpha_{S_m})(f)\|_2<\eps.
$$
We assert that if $\eps<\delta^2/4$ then there exist $g\in G$ and
$c\in\T$ such that
\begin{equation} \label{e5.1}
\|w_n-cw_mv_g\|_2<(2\eps)^{1/2},
\end{equation}
where $v_g$, $g\in G$, are the canonical generators of $S_T$.
Indeed, let $v=w_m^*w_n$, $v=\sum_ga_gv_g$, $a_g\in\C$. If
$E\colon R_T\to C_T$ is the trace-preserving conditional
expectation, then for arbitrary $x\in C_T$ we have
$E(vxv^*)=\sum_g|a_g|^2\alpha_g(x)$, whence
\begin{eqnarray*}
\eps^2
 &>&\|(\Ad w_n\circ\alpha_{S_n})(f)
      -(\Ad w_m\circ\alpha_{S_m})(f)\|_2^2\\
 &=&2(1-\RR\tau(v\alpha_{S_n}(f)v^*\alpha_{S_m}(f^*)))\\
 &=&2(1-\RR\tau(E(v\alpha_{S_n}(f)v^*)\alpha_{S_m}(f^*)))\\
 &=&2\sum_g|a_g|^2(1-\RR\tau((\alpha_g\circ\alpha_{S_n})(f)
      \alpha_{S_m}(f^*))).
\end{eqnarray*}
Set $Y=\{g\in G\ |\
1-\RR\tau((\alpha_g\circ\alpha_{S_n})(f)\alpha_{S_m}(f^*))<\eps/2\}$.
If $g\in Y$ then
$\|(\alpha_g\circ\alpha_{S_n})(f)-\alpha_{S_m}(f)\|_2<\eps^{1/2}$.
Thus if $g_1\ne g_2$ both lie in $Y$ then
$\|(\alpha_{g_1}\circ\alpha_{S_n})(f)
 -(\alpha_{g_2}\circ\alpha_{S_n})(f)\|_2<2\eps^{1/2}$.
Since
$\alpha_g\circ\alpha_{S_n}=\alpha_{S_n}\circ\alpha_{\beta(g)}$ for
some automorphism $\beta$, we get a contradiction if
$2\eps^{1/2}<\delta$. Hence the set $Y$ consists of at most one
point. On the other hand, we have
$$
\eps^2>2\sum_g|a_g|^2(1
 -\RR\tau((\alpha_g\circ\alpha_{S_n})(f)\alpha_{S_m}(f^*)))
\ge\eps\sum_{g\notin Y}|a_g|^2,
$$
so that $\sum_{g\notin Y}|a_g|^2<\eps$. It follows that $Y$ is
non-empty. Hence it consists precisely of one point~$g$, and
$|a_g|^2>1-\eps$. Then with $c=a_g/|a_g|$ we have
$$
\|w_n-cw_mv_g\|_2^2=\|v-cv_g\|_2^2
 =\sum_{h\ne g}|a_h|^2+|a_g-c|^2
 =2(1-|a_g|)<2\eps.
$$

It follows that passing to a subsequence we may suppose that for
each $n\ge2$ there exist $g_n\in G$ and $c_n\in\T$ such that
$$
\|c_nw_nv_{g_n}-w_{n-1}\|_2<{1\over2^n}.
$$
Then replacing $w_n$ by $c_n\ldots c_2w_nv_{g_n+\ldots+g_2}$ and
$S_n$ by $T_{-g_2-\ldots-g_n}S_n$ we still have $\alpha_n=\Ad
w_n\circ\alpha_{S_n}$, but now the sequence $\{w_n\}_n$ converges
strongly to a unitary $w\in S_T$. Then $\alpha(C_T)=(\Ad w)(C_T)$.
\endverif

Note that the group $\AutTT$ consisting of all automorphisms
$\gamma\in\Aut(R_T,C_T)$ such that $\gamma(S_T)$ and $S_T$ are
inner conjugate, is never closed. Indeed, let $c\in Z^1(\RT,\T)$
be a $\T$-valued 1-cocycle on the orbit equivalence relation $\RT$
defined by $T$, and $\sigma_c\in\Aut(R_T,C_T)$ the corresponding
automorphism~\cite{FM2}. Then $\sigma_c(S_T)$ and $S_T$ are inner
conjugate if and only if $c$ is cohomologous to the cocycle
$c_\chi$, $c_\chi(x,T_gx)=\<\chi,g\>$, for some $\chi\in\hat
G$~\cite{Packer1} (this result was proved in~\cite{Packer1} for
actions with purely discrete spectrum, but with minor changes the
proof works for arbitrary ergodic actions; in our weakly mixing
case using Theorem~\ref{4.2} and the fact that if
$\gamma\in\AutTT$ then $\Ad u\circ\gamma\in\AutT$ for some unitary
$u$, it is easy to obtain a more precise result: the group
$\AutTT$ consists of automorphisms of the form
$\sigma_c\circ\alpha_S$, where $c$ is a cocycle cohomologous to
$c_\chi$ and $S\in I(T)[T]$, where $[T]$ is the full group
generated by $T_g$, $g\in G$). Since the equivalence relation is
hyperfinite, any cocycle can be approximated by coboundaries, so
all automorphisms $\sigma_c$ are in the closure of $\AutTT$. On
the other hand, there always exist cocycles which are not
cohomologous to cocycles $c_\chi$, because otherwise $Z^1(\RT,\T)$
would be a continuous isomorphic image of the group $\hat G\times
I(X,\T)$, where $I(X,\T)$ is the factor of the unitary group of
$L^\infty(X)$ by the scalars (note that since the action is weakly
mixing, $c_\chi$ is not a coboundary for $\chi\in\hat
G\backslash\{e\}$), hence $Z^1(\RT,\T)$ would be topologically
isomorphic to $\hat G\times I(X,\T)$, which would imply that the
group of coboundaries is closed.

\smallskip

If the action is rigid, it is still possible that $\AutT$ is
closed. However, as the following result shows the group
$\Int(S_T)$ consisting of inner automorphisms of $R_T$ defined by
unitaries in $S_T$ is not closed in this case, which may indicate
that we should consider systems satisfying stronger mixing
properties than weak mixing. Note that if an action is mixing then
it is not rigid.

\begin{proposition} \label{4.11}
The following conditions are equivalent:
\newline
{\rm(i)} the action $T$ is rigid;
\newline
{\rm(ii)} there exist non-trivial central sequences in $S_T$;
\newline
{\rm(iii)} the subgroup $\Int(S_T)$ of $\Aut(R_T)$ is not closed.
\end{proposition}

\noindent{\it Proof.} The equivalence of (ii) and (iii) is
well-known~\cite{Connes}. The implication (i)$\Rightarrow$(ii) is
obvious. Suppose that the action is not rigid. Let $\{u_n\}_n$ be
a central sequence  of unitaries in $S_T$. For fixed $n$ apply
(\ref{e5.1}) to $w_n=u_n$, $w_m=1$, $S_n=S_m=\id$. Then we
conclude that there exist $c_n\in\T$ and $g_n\in G$ such that
$\|u_n-c_nv_{g_n}\|_2\to0$ as $n\to\infty$. The sequence
$\{v_{g_n}\}_n$ is central, which is equivalent to the strong
convergence $u_{g_n}\to1$. Since the action is not rigid, this
implies that eventually $g_n=0$, so the central sequence
$\{u_n\}_n$ is trivial. Thus (ii) implies (i).
\endverif

The following corollary is not surprising in view of
Proposition~\ref{3.1} but is worth mentioning.

\begin{corollary}
There exist weakly mixing transformations $T^{(1)}$ and $T^{(2)}$
such that the singular masas $S_{T^{(1)}}$ and $S_{T^{(1)}}$ are
not conjugate but their Pukanszky invariants coincide.
\end{corollary}

\noindent{\it Proof.} The class of weakly mixing
measure-preserving transformations with simple spectrum, i.e. of
spectral multiplicity one, contains both rigid and non-rigid
transformations (e.g. certain Gauss systems are rigid and have
simple spectrum~\cite[Chapter~14]{CFS}, while Ornstein's rank-one
transformations are mixing~\cite[Chapter~16]{Nadkarni}). Since
rigidity is a conjugacy invariant by Proposition~\ref{4.11}, there
exist transformations $T^{(1)}$ and $T^{(2)}$ such that
$S_{T^{(1)}}$ and $S_{T^{(1)}}$ are not conjugate, while
$P(S_{T^{(1)}})=P(S_{T^{(2)}})=\{1\}$.
\endveriff

\bigskip\bigskip

\section{Entropy} \label{5}

A weak form of our conjecture would be to say that conjugacy of
masas $S_T$ for actions of an abelian group $G$ implies
coincidence of the entropies. In this form the conjecture may hold
without any assumptions on the spectrum, since systems with purely
discrete spectrum have zero entropy. The main result of this
section is a step towards the solution of this weaker problem.
While in the previous section we proved that if the conjecture is
false then the isomorphism $\gamma\colon R_{T^{(1)}}\to
R_{T^{(2)}}$ for non-isomorphic systems sends $C_{T^{(1)}}$ far
from $C_{T^{(2)}}$, in this section we shall prove that if the
entropies are distinct, the images $\gamma(v_g)$ of the canonical
generators of $S_{T^{(1)}}$ can not coincide with the generators
of $S_{T^{(2)}}$ even on small projections.

We shall consider only the case $G=\Z$, since the theory of
non-commutative entropy is not well-developed for actions of
general abelian (or amenable) groups, though in fact the result is
true for arbitrary abelian~$G$.

\begin{theorem} \label{5.1}
Let $T^{(i)}\in\Aut(X_i,\mu_i)$ be a measure-preserving
transformation, $i=1,2$. Denote by $v_i$ the canonical generator
of $S_{T^{(i)}}$. Suppose there exists an isomorphism
$\gamma\colon R_{T^{(1)}}\to R_{T^{(2)}}$ such that
$\gamma(S_{T^{(1)}})=S_{T^{(2)}}$, and the unitary
$\gamma(v_1)v_2^*$ has an eigenvalue. Then
$h(T^{(1)})=h(T^{(2)})$.
\end{theorem}

The result will follow from

\begin{proposition} \label{5.2}
Let $T\in\Aut(X,\mu)$ be a measure-preserving transformation,
$v\in S_T$ the canonical generator. Then for any non-zero
projection $p\in S_T$ we have $H(\Ad v|_{pR_Tp})=h(T)$, where
$H(\Ad v|_{pR_Tp})$ is the entropy of Connes and
St{\o}rmer~\cite{CS} of the inner automorphism $\Ad v|_{pR_Tp}$
computed with respect to the normalized trace
$\tau_p=\tau(p)^{-1}\tau|_{pR_Tp}$.
\end{proposition}

\noindent{\it Proof of Theorem~\ref{5.1}.} By assumption, there
exists $\theta\in\T$ such that the spectral projection~$p$ of the
unitary $\gamma(v_1)v_2^*$ corresponding to the set $\{\theta\}$
is non-zero. Then $\gamma(v_1)p=\theta v_2p$. By
Proposition~\ref{5.2} we get
$$
h(T_1)=H(\Ad v_1|_{\gamma^{-1}(p)R_{T^{(1)}}\gamma^{-1}(p)})
=H(\Ad\gamma(v_1)|_{pR_{T^{(2)}}p})) =H(\Ad
v_2|_{pR_{T^{(2)}}p}))=h(T_2).
$$
\endveriff

To prove Proposition~\ref{5.2} consider a more general situation
when we are given a finite injective von Neumann algebra $M$ with
a fixed normal faithful trace $\tau$ and a $\tau$-preserving
automorphism~$\alpha$. For each projection $p$ in the fixed point
algebra $M^\alpha$ we set
$$
\tau_\alpha(p)=\tau(p)H(\alpha|{pMp}).
$$

\begin{proposition} \label{5.3}
The mapping $p\mapsto\tau_\alpha(p)$ extends uniquely to a normal
(possibly infinite) trace $\tau_\alpha$ on~$M^\alpha$, which is
invariant with respect to all $\tau$-preserving automorphisms in
$\Aut(M,M^\alpha)$ commuting with~$\alpha$.
\end{proposition}

\noindent{\it Proof.} To prove that the mapping extends to a
normal trace it is enough to check that the following three
properties are satisfied: $\tau_\alpha(upu^*)=\tau_\alpha(p)$ for
any unitary $u$ in~$M^\alpha$, if $p_n\nearrow p$ then
$\tau_\alpha(p_n)\nearrow\tau_\alpha(p)$, the mapping
$p\mapsto\tau_\alpha(p)$ is finitely additive.

The first property is a particular case of the last statement of
the proposition. If $\beta\in\Aut(M,M^\alpha)$ commutes with
$\alpha$ and preserves the trace $\tau$, then it defines an
isomorphism of the systems $(pMp,\tau_p,\alpha)$ and
$(\beta(p)M\beta(p),\tau_{\beta(p)},\alpha)$, so their entropies
coincide.

The second property follows from the well-known continuity
properties of entropy:
$$
\tau_\alpha(p_n)=\tau(p_n)H(\alpha|{p_nMp_n})
=\tau(p)H(\alpha|_{p_nMp_n+\C(p-p_n)})\nearrow
\tau(p)H(\alpha|_{pMp})=\tau_\alpha(p).
$$

To prove the third one consider a finite family $\{p_i\}^n_{i=1}$
of mutually orthogonal projections in~$M^\alpha$ and set
$p=\sum_ip_i$. Let
$$
B=p_1Mp_1+\ldots+p_nMp_n.
$$
By affinity of entropy,
$$
H(\alpha|_B)=\sum_i{\tau(p_i)\over\tau(p)}H(\alpha|_{p_iMp_i})
=\tau(p)^{-1}\sum_i\tau_\alpha(p_i).
$$
So in order to prove finite additivity it is enough to prove that
$H(\alpha|_{pMp})=H(\alpha|_B)$. The trace-preserving conditional
expectation $E\colon pMp\to B$ has the form
$$
E(x)=p_1xp_1+\ldots+p_nxp_n.
$$
It commutes with $\alpha$ and is of finite index, $E(x)\ge {1\over
n}x$ for $x\in pMp$, $x\ge0$. Indeed, if we consider $pMp$ acting
on some Hilbert space, then for a vector $\xi$ we set
$\xi_i=p_i\xi$ and get
\begin{eqnarray*}
(x\xi,\xi)
 &=&\sum_{i,j}(x^{1/2}\xi_i,x^{1/2}\xi_j)
      \le\sum_{i,j}\|x^{1/2}\xi_i\|\cdot\|x^{1/2}\xi_j\|
      =\left(\sum_i\|x^{1/2}\xi_i\|\right)^2\\
 &\le&n\sum_i\|x^{1/2}\xi_i\|^2=n\sum_i(xp_i\xi,p_i\xi)=n(E(x)\xi,\xi).
\end{eqnarray*}
By \cite[Corollary 2]{NS}, we conclude that
$H(\alpha|_{pMp})=H(\alpha|_B)$.
\endverif

\noindent{\it Proof of Proposition~\ref{5.2}.} Consider the weight
$\tau_{\Ad v}$ on $S_T$ corresponding to the automorphism $\Ad v$
of $R_T$. Then we have to prove that $\tau_{\Ad
v}=h(T)\tau|_{S_T}$. By Proposition~\ref{5.3}, the weight
$\tau_{\Ad v}$ is invariant under the dual action. Since this
action is ergodic on $S_T$, $\tau_{\Ad v}$ is a scalar multiple of
$\tau|_{S_T}$, $\tau_{\Ad v}=c\cdot\tau|_{S_T}$ for some
$c\in[0,+\infty]$. By definition of $\tau_{\Ad v}$ we have
$c=H(\Ad v)$. But by~\cite{GN,Vo3}, $H(\Ad v)=h(T)$, and the proof
is complete.
\endverif

The definition of the weight above leads to the
following interesting problem in entropy theory. Let $A$ be an
abelian subalgebra of a finite algebra $M$. For each unitary $u\in
A$ consider the weight $\tau_u$ on $A$, which is the restriction
of the weight $\tau_{\Ad u}$ to~$A$.

\medskip\noindent
{\bf Problem.} Find the connection between $\tau_u$ and
$\tau_{\phi(u)}$, where $\phi$ is a Borel mapping from $\T$ onto
itself.

\medskip

Voiculescu's approach to entropy using norm of
commutators~\cite{Vo1, Vo2} suggests that such a connection exists
at least when $\phi$ is smooth. More interesting is the case when
$u$ is a Haar unitary and $\phi$ is an invertible transformation
preserving Lebesgue measure, so that $\phi(u)$ is again Haar and
generates the same algebra. Note also some resemblance of this
problem to the computation of entropy of Bogoliubov
automorphisms~\cite{SV,N}. However, the correspondence
$u\mapsto\tau_u$ does not have nice continuity properties which
makes the problem more difficult.

\medskip

Finally note that the problems studied in the paper can also be considered
for topological dynamical systems and C$^*$-crossed products. In this
setting isomorphism of crossed products already implies that the systems
have a non-trivial relationship. For example, for minimal homeomorphisms
of Cantor sets the crossed products are isomorphic if and only if the
systems are strongly orbit equivalent~\cite{GPS}. Since rotations are
the only measure- and orientation-preserving homeomorphisms of the circle,
if $\gamma$ is an isomorphism of $C(X_1)\rtimes\Z$ onto
$C(X_2)\rtimes\Z$ which maps $C^*(v_1)$ onto $C^*(v_2)$ then 
$\gamma(v_1)=\theta v_2^{\pm1}$, so the homeomorphisms have the same topological
entropy.

\bigskip\bigskip

\begin{flushleft}
University of Oslo\\
Blindern, P.O. Box 1053\\
0316 Oslo, Norway\\
{\it e-mail}: neshveyev@hotmail.com, erlings@math.uio.no
\end{flushleft}

\end{document}